\input amstex
\documentstyle{amsppt}
\magnification=\magstep1

\pageheight{9.0truein}
\pagewidth{6.7truein}

\NoBlackBoxes
\TagsAsMath

\input xy \xyoption{matrix} \xyoption{arrow} \xyoption{curve}

\define\seq{\mathrel{\widehat{=}}}
\define\la{{\Lambda}}
\define\lamod{\Lambda\text{-}\roman{mod}}
\define \len{\operatorname{length}} 
\define\bbA{{\Bbb A}}
\define\PP{{\Bbb P}}
\define\SS{{\Bbb S}}

\define\NN{{\Bbb N}}
\define\RR{{\Bbb R}}
\define\aut{\operatorname{Aut}}

\define\GL{\operatorname{GL}}

\define\Ker{\operatorname{Ker}}
\define\Img{\operatorname{Im}}

\define\autlap{\operatorname{Aut}_\la(P)}
\define\End{\operatorname{End}}

\define\A{{\Cal A}}
\define\C{{\Cal C}}




\define\R{{\Cal R}}

\define\frakt{{\frak t}}

\define\e{\bold {e}}



\define\lay#1{\operatorname{\bold{Mod}}(#1)}

\define\grasstd{\operatorname{\frak{Grass}}^T_d}
\define\Gr{\operatorname{\frak{Gr}}}

\define\ggrasstd{\operatorname{Gr-\frak{Grass}}^T_d}
\define\gautlap{\operatorname{Gr-\aut_\Lambda(P)}}

\define\grassS{\operatorname{\frak{Grass}}(\sigma)}
\define\ggrassS{\operatorname{Gr-\frak{Grass}}(\sigma)}
\define\grassSS{\operatorname{\frak{Grass}}(\SS)}
\define\ggrassSS{\operatorname{Gr-\frak{Grass}}(\SS)}
\define\grass#1{\operatorname{\frak{Grass}}(#1)}
\define\ggrass#1{\operatorname{Gr-\frak{Grass}}(#1)}

\define\Schu{\operatorname{Schu}}

\define\GeomIV{{\bf 1}}
\define\Hil{{\bf 2}}
\define\GeomI{{\bf 3}}
\define\grassI{{\bf 4}}
\define\grassII{{\bf 5}}
\define\grassIII{{\bf 6}}
\define\Ki{{\bf 7}}
\define\New{{\bf 8}}

\topmatter

\title Moduli spaces of graded representations of finite dimensional
algebras 
\endtitle

\rightheadtext{Graded representations}

\author E. Babson, B. Huisgen-Zimmermann, and R. Thomas
\endauthor

\address Department of Mathematics,
University of Washington, Seattle, Washington 98195 \endaddress

\email babson\@math.washington.edu \endemail

\address Department of Mathematics, University of California,
Santa Barbara, California 93106 \endaddress

\email birge\@math.ucsb.edu \endemail
 
\address Department of Mathematics,
University of Washington, Seattle, Washington 98195 \endaddress

\email thomas\@math.washington.edu \endemail



\thanks The authors were partially supported by grants
from the National Science Foundation.  Moreover, the second-named author
would like to thank the University of Washington for the hospitality
extended to her while this work was in progress. 
\endthanks

\abstract Let $\la$ be a basic finite dimensional algebra over an
algebraically closed field, presented as a path algebra modulo
relations; further, assume that $\la$ is graded by lengths of paths. The
paper addresses the classifiability, via moduli spaces, of classes of
graded $\la$-modules with fixed dimension $d$ and fixed top $T$. It is
shown that such moduli spaces exist far more frequently than they do for
ungraded modules. In the local case (i.e., when $T$ is
simple), the graded $d$-dimensional $\la$-modules with top $T$ always
possess a fine moduli space which classifies these modules up to
graded-isomorphism; moreover, this moduli space is a projective variety
with a distinguished affine cover that can be constructed from quiver and
relations of $\la$. When $T$ is not simple, existence of a
coarse moduli space for the graded $d$-dimensional $\la$-modules with top
$T$ forces these modules to be direct sums of local modules; under the
latter condition, a finite collection of isomorphism invariants of the
modules in question yields a partition into subclasses, each of which has
a fine moduli space (again projective) parametrizing the
corresponding graded-isomorphism classes.
\endabstract

\endtopmatter

\document

\head 1. Introduction
\endhead

Let $\la$ be a finite dimensional algebra with radical $J$ over an
algebraically closed field $K$, and fix a finite
dimensional semisimple (left)
$\la$-module $T$ together with a positive integer $d$.  In
\cite{\grassI}, the second author explored the existence and structure
of moduli spaces classifying, up to isomorphism, those
$d$-dimensional (left) representations $M$ of $\la$ whose tops 
$M/JM$ equal $T$ under identification of isomorphic semisimple
modules.  The vehicle for tackling this classification problem is a
projective variety, $\grasstd$, parametrizing the $d$-dimensional
modules with top $T$; see Section 2. 

The general goal driving such investigations is to demonstrate
that, even over a wild algebra $\la$, major portions of the
representation theory may behave tamely, being accessible to
classification in a quite stringent sense.  The idea of {\it moduli\/}
goes back to Riemann's 1857 classification of nonsingular projective
curves of fixed genus in terms of continuous structure determining
invariants; it was made precise by Mumford in the 1960's.  In rough
terms, adapted to our representation-theoretic context, it amounts to
the following:  First one introduces a concept of {\it family\/}, which
pins down what it means that a collection of finite dimensional
representations be parametrized ``continuously" by the points of a
variety $X$.  Given this prerequisite, a fine or coarse moduli space for
a class of finite dimensional representations of $\la$ is a variety that
continuously and bijectively parametrizes the isomorphism classes of
the considered representations in a fashion satisfying a certain
{\it coarse\/} or {\it fine\/} universal property.  In slightly more
precise language, a {\it fine moduli space\/}  --  the crucial concept
in this paper  -- is the parametrizing variety of a
distinguished family satisfying the postulate that any other family of
representations recruited from the given class be uniquely ``induced"
by the distinguished one.  Existence provided, both fine and coarse
moduli spaces are unique up to canonical isomorphism due to the
pertinent universality conditions.      

Here we focus on a {\it graded\/} basic finite dimensional
algebra
$\la$ and address two problems closely related to the one
mentioned at the outset:  (1)  That of deciding classifiability of the
{\it graded\/} $d$-dimensional left $\la$-modules $M$ with fixed top
$T$, up to graded-isomorphism; and, more restrictively,
classifiability of those graded candidates $M$ which have fixed
{\it radical layering\/}
$\SS(M) = \bigl(J^l M/J^{l+1} M \bigr)_{l \ge 0}$.  In either case,
``classifiability" stands for existence of a fine or coarse moduli
space.  (2) In case existence of a moduli space is secured, the
problem of determining the structure of this space and of
constructing a universal family for the  considered class of
representations.   

Our base field $K$ being algebraically closed, we may assume without
loss of generality that $\la$ is a path algebra modulo relations,
meaning that
$\la = KQ/I$ for a quiver
$Q$ and an admissible ideal $I$ in the path algebra $KQ$.  In fact, we
specialize to the situation where $\la$ is {\it graded by lengths of
paths\/}, meaning that $I$ is a homogeneous ideal with respect to the
natural grading of $KQ$ through path lengths.  We will start by showing
that the set of those points in the mentioned variety $\grasstd$, which
correspond to the graded
$d$-dimensional modules with top $T$ that are generated in degree zero,
form a closed  --  and hence projective  --  subvariety of $\grasstd$, 
denoted by $\ggrasstd$.  More strongly, we will verify the
following:  Suppose $\SS = (\SS_0, \SS_1, \dots, \SS_L)$
is a sequence of semisimple modules with $\SS_0 = T$, where $J^{L+1} =
0$ and the dimensions of the $\SS_l$ add up to $d$.  Then the following
subset $\ggrassSS$ of
$\grasstd$ is closed:  Namely, the set of those points in 
$\ggrasstd$ which correspond to the graded modules $M$ with 
$\SS(M) = \SS$.  In alternate terms, $\ggrassSS$ is a projective
variety parametrizing the
$d$-dimensional graded modules generated in degree zero with 
radical layering
$\SS$.  Closedness of these subvarieties entails, in particular, that
each
$\ggrassSS$ is a union of irreducible components of $\ggrasstd$.  This
is the first crucial difference between the graded and ungraded
settings.  Indeed, by contrast, the subvariety 
$\grassSS$ consisting of {\it all\/} points in
$\grasstd$ corresponding to (not necessarily graded) modules with
radical layering
$\SS$ fails to be closed in $\grasstd$ in general.

Naturally, the {\it graded\/} $d$-dimensional representations with top
$T$ possess a fine/{\allowlinebreak}coarse moduli space whenever all
$d$-dimensional representations with top $T$ do.  On the other hand, not
too surprisingly, existence of a moduli space is a far more frequent event
in the graded than in the ungraded situation, as a grading accounts
for increased rigidity.  What is surprising is the extent of this
discrepancy:   For instance, given a simple module
$T$ with projective cover $P$, the  $d$-dimensional top-$T$ modules have
a fine (equivalently, a coarse) moduli space classifying them up to
isomorphism precisely when every submodule $C$ of $JP$ having
codimension $d$ in $P$ is invariant under endomorphisms of $P$; the
latter requirement imposes strong restrictions on the underlying triple 
$(\la, T, d)$; see \cite{\grassI\rm, Corollary 4.5}.  However, when one
narrows one's view to graded representations under graded-isomorphism,
existence of a fine moduli space is automatic for a simple top $T$:  

\proclaim{Theorem A}  If $\la$ is path-length-graded and $T$ a
simple left $\la$-module, then, for any positive integer $d$,
the graded $d$-dimensional $\la$-modules with top $T$ possess a fine
moduli space, classifying their graded-isomorphism classes. 
This moduli space equals $\ggrasstd$.
\endproclaim  

Calling a module {\it local\/} if it has a simple top, we will more generally
prove the following:

\proclaim{Theorem B}  Suppose that $\la$ is path-length-graded, $T \in
\lamod$ any semisimple $\la$-module, and $\SS$ a sequence of semisimple
$\la$-modules as above.  Moreover, let $\C(T)$ {\rm{(}}resp\.
$\C(\SS)${\rm{)}} be the class of all graded
$d$-dimensional $\la$-modules with top $T$ {\rm{(}}resp\. with radical
layering $\SS${\rm{)}}.  Then the following are true:

$\bullet$ If there is a coarse
moduli space classifying the graded-isomorphism classes in $\C(T)$
{\rm{(}}resp\. $\C(\SS)${\rm{)}}, then every object in $\C(T)$
{\rm{(}}resp\. $\C(\SS)${\rm{)}} is a direct sum of local modules.

$\bullet$  Conversely, if $\C(T)$ {\rm{(}}resp\. $\C(\SS)${\rm{)}}
consists of direct sums of local modules, then $\C(T)$ {\rm{(}}resp\.
$\C(\SS)${\rm{)}} can be partitioned into finitely many subclasses,
each of which has a fine moduli space.

All moduli spaces arising in the latter case are
projective.   
\endproclaim 

In parallel with the ungraded situation, each of the varieties
$\ggrasstd$ possesses a distinguished affine cover, accessible from
quiver and relations of $\la$, which provides the key to analyses of
concrete examples.     

This leaves the question of which projective varieties
occur among the irreducible components of fine moduli spaces for
graded modules with fixed dimension and top.  We use examples of
Hille in \cite{\Hil}, which are in turn based on a construction
technique introduced by the second author in \cite{\GeomI}, to
show that every irreducible projective variety arises as an irreducible
component of such a space.

Our approach to moduli problems for representations is fundamentally
different from that of King in \cite{\Ki}, where the  targeted modules 
are those that are semistable with respect to a given
additive function $\Theta:  \operatorname{K}_0 (\lamod) \rightarrow
\RR$.  King's definition of semistability allows for the adaptation of
techniques developed by Mumford with the aim of classifying vector
bundles.  On one hand, in King's approach (coarse) moduli spaces for
$\Theta$-semistable representations are guaranteed to exist.  On the
other hand, in general these classes of modules are hard to assess in
size and to describe in more manageable terms, while
their classification through moduli spaces is a priori only up to an
equivalence relation considerably coarser than isomorphism.

Concerning strategy:  Evidently, every
local graded module is generated in a single degree, which, for
purposes of classification, we may assume to be zero.   As for the
general case, we will show that classifiability up to
graded-isomorphism (through a moduli space) of the
graded $d$-dimensional modules with fixed top $T$, generated in mixed
degrees, forces these graded objects to be direct sums of
local graded submodules.  We are thus led back to a situation in which
restriction to graded modules generated in degree zero is harmless.  The proof of
this reduction step requires an extra layer of technicalities likely to
obscure the underlying ideas; we will therefore defer it to an
appendix (Section 6).  In Sections 2--5, we will only consider graded modules
generated in degree zero.  

In Section 2, we will provide prerequisites; in particular, we will
introduce the varieties
$\ggrasstd$ and
$\ggrassSS$ and verify their projectivity.  In Section 3, we will
prepare for proofs of the main results by introducing the mentioned
affine cover of the variety $\ggrasstd$ and by constructing a pivotal
family of graded modules with top $T$; this family will turn out to be
the universal one (see Section 2 for a definition) in case a fine moduli
space exists. Section 4 contains proofs of upgraded versions of
Theorems A and B, the latter restricted to graded modules generated in
degree zero.  Section 5 is devoted to examples.  The appendix, finally,
will remove the restriction concerning degree-zero generating sets from
the results for nonlocal modules. 
\bigskip              

\head 2. Further terminology and Background \endhead

 We will be fairly complete in setting up our conventions, even
fairly standard ones, for the convenience of the reader whose
expertise lies at the periphery of the subject. 
  
Let $\la$ be a basic 
finite dimensional algebra over an algebraically closed field $K$. 
Without loss of generality, we assume $\la$
to be a path algebra modulo relations, that is, $\la = KQ/I$ for a
quiver $Q$ and an admissible ideal $I$ in the path algebra
$KQ$. 

{\it Gradings.\/} Throughout, we suppose $\la$ to be graded in terms of
path lengths, meaning that
$I$ is homogeneous with respect to the length-grading of
$KQ$.  Denoting by $J$ the Jacobson radical of $\la$, we let $L$ be
maximal with $J^L \ne 0$.  Then the grading of $\la$ takes on the form
$\la = \bigoplus_{0 \le l \le L} \la _l$, where $\la_l \cong
J^l/J^{l+1}$ is the homogeneous component of degree $l$ of
$\la$.  The vertices
$e_1, \dots, e_n$ of
$Q$ will be identified with the primitive idempotents of $\la$
corresponding to the paths of length zero, that is, the $e_i$ will also
stand for the $I$-residues of the paths of length $0$ in $\la_0$.
The factor modules $S_i =
\la e_i /J e_i$ then form
an irredundant set of representatives for the simple left
$\la$-modules; unless we explicitly state otherwise, we consider the
$S_i$  --  and hence all semisimple modules -- as homogeneous modules in
degree $0$, systematically identifying isomorphic semisimple modules. 
Clearly, the grading of any indecomposable projective module
$\la e_i$ inherited from that of
$\la$ yields a graded local module which is generated in degree zero. 
Whenever $P = \bigoplus_{1
\le i \le n} (\la e_i)^{t_i}$, we let $P = \bigoplus_{0 \le l \le L}
P_l$ be the resulting decomposition into homogeneous subspaces.  Given
two graded modules $M$, $M'$, we call a morphism $f: M \rightarrow M'$
{\it homogeneous of degree\/} $s$ in case $f(M_l) \subseteq M'_{l +
s}$ for all $l$; the attribute ``homogeneous" by itself stands for
``homogeneous of degree zero".  Whenever there is an isomorphism $M
\rightarrow M'$ which is homogeneous of some degree $s$, we call
$M$ and $M'$ {\it graded-isomorphic\/}; so, in particular, two graded
modules generated in degree zero are graded-isomorphic if and only
if they are isomorphic by way of a homogeneous map.   

{\it Paths in $\la$ and top elements of modules.\/} We
will observe the following conventions:  The product $pq$ of two paths
$p$ and $q$ in $KQ$  stands for ``first
$q$, then $p$"; in particular, $pq$ is zero unless the end point of $q$
coincides with the starting point of $p$.  In this spirit, we call a
path
$p_1$ a {\it right subpath\/} of $p$ if $p = p_2 p_1$ for some path
$p_2$.  A {\it path in
$\la$\/} will be an element of the form $p +I$, where
$p$ is a path in $K Q \setminus I$, whence $\la_l$ consists of the
$K$-linear combinations of paths of length $l$ in $\la$; in light of
homogeneity of the ideal $I$, the length of a path in $\la$ is an
unambiguous quantity.  We will suppress the residue notation, provided
there is no risk of confusion, and gloss over the
distinction between the left
$\la$-structure of a module $M \in \lamod$ and its induced
$KQ$-structure.  A homogeneous element
$x$ of a graded module $M$ generated in degree zero will be
called a {\it top element\/} of
$M$ if $x \in M_0 \setminus JM$ and $x$ is normed by some $e_i$, meaning
that $x = e_i x$.  If $M$ is ungraded, we waive the homogeneity
condition imposed on $x$.  Any sequence
$x_1, \dots, x_m$ of top elements of
$M$ generating
$M$ and linearly independent modulo $JM$ will be called a {\it full
sequence of top elements of $M$\/}. 

{\it Background on moduli problems.\/} We refer to
\cite{\New}, but recall the definition of a fine moduli space in the
context of representations.  Our concept of a
family of
$\la$-modules is that introduced by King in
\cite{\Ki}:  Namely, a {\it family of $d$-dimensional $\la$-modules 
parametrized by an algebraic variety
$X$\/} is a pair $(\Delta, \delta)$, where $\Delta$ is a (geometric)
vector bundle of rank $d$ over $X$ and $\delta: \la \rightarrow
\End(\Delta)$ a $K$-algebra homomorphism.  We are interested in {\it
families of graded $d$-dimensional representations of $\la$\/}, which
means that the left
$\la$-multiplication induced by $\delta$ on each fibre $\Delta_x$
of $\Delta$ yields a graded module.  In addition, we want the
$\Delta_x$ to be ``continuously graded" as $x$ traces $X$, meaning that
all fibres are generated in the same sets of degrees (see Section 6). 
Primarily (Sections 2--5), we will be interested in families of
$\la$-modules all of which are generated in degree zero.   Our
notion of {\it equivalence of families of graded modules\/}
parametrized by some variety 
$X$ is the coarsest possible to separate
graded-isomorphism classes:  Namely,
$(\Delta^1, \delta^1) \sim (\Delta^2, \delta^2)$ precisely when, for
each $x \in X$, the fibre of $\Delta^1$ above $x$ is graded-isomorphic
as a left $\la$-module to the fibre of $\Delta^2$ above $x$.  As is
common, given a family
$(\Delta, \delta)$ parametrized by
$X$ and a morphism $\tau: Y
\rightarrow X$ of varieties, the {\it induced family\/} $\tau^*(\Delta,
\delta)$ over $Y$ is the pullback of $(\Delta, \delta)$ along $\tau$.
The family $\tau^*(\Delta, \delta)$ is again a family of graded
modules if $(\Delta, \delta)$ is.  In this setup, a variety
$Z$ is a {\it fine moduli space\/} for the (families of)
graded $d$-dimensional modules with top $T$ if there exists a family
$(\Gamma, \gamma)$ of such modules parametrized by $Z$,
which has the property that an arbitrary family parametrized by some
variety $X$ is equivalent to a family
$\tau^*(\Gamma, \gamma)$ induced via a unique morphism $\tau: X
\rightarrow Z$.  Accordingly, $\Gamma$ is then called the {\it
universal\/} family.  In particular, the
requirements on
$\Gamma$ entail that every
$d$-dimensional module with top $T$ be isomorphic to precisely one
fibre of the bundle
$\Gamma$.  More commonly (but equivalently), a fine moduli space is
defined in terms of representability of the contravariant functor from
the category of varieties to the category of sets, assigning to any
variety $X$ the set of equivalence classes of families parametrized
by $X$.  For details, as well as for the concept of  a {\it coarse
moduli space\/}, see
\cite{\New\rm, pp\. 23, 24}. 
\smallskip

{\it Ungraded Grassmannians of modules with top $T$.\/} \ We fix a
natural number
$d$, a semisimple module
$T$, say $T = \bigoplus_{1 \le i \le n} S_i^{t_i}$, 
and denote by $P = \bigoplus_{1 \le i \le n} (\la
e_i)^{t_i}$ its projective cover.  The module $P$ comes equipped with
the obvious grading under which it is obviously generated in degree
zero.  It will be convenient to write $P$ in the form  
$$P = \bigoplus_{1 \le r \le \frakt} \la z_r,$$ 
where $\frakt = \sum_i t_i$ and
$z_1, \dots, z_t$ is a sequence of top elements (all of which are
homogeneous of degree zero); in particular, $ z_r = e(r)z_r$ for
suitable idempotents $e(r)$ in
$\{e_1, \dots, e_n\}$.  A natural choice of
such top elements consists of the  primitive idempotents $e_i$
themselves, each with multiplicity
$t_i$, distinguished by their `slots' in the above decomposition of
$P$; to distinguish these slots, we will write $\e(r)$ for the element
of the direct sum $\bigoplus_{1 \le i \le n} (\la e_i)^{t_i}$ which has
entries zero outside the slot labeled $r$ and carries $e(r)$ in the
$r$-th slot.   

The ungraded Grassmannian of $d$-dimensional left $\la$-modules
with top $T$ was defined in \cite{\grassI} as follows: If
$\Gr(d',JP)$ is the classical Grassmannian of $d'$-di\-men\-sion\-al
subspaces of the $K$-space $JP$, where $d' = \dim_K P
- d$, then
$$\grasstd = \{C \in \Gr(d',JP) \mid C \text{\ is a\ }
\la\text{-submodule of\ } JP \},$$ 
a closed subvariety of $\Gr(d',JP)$ and consequently
projective.  It is accompanied by a natural surjection
$$\grasstd \longrightarrow \{\text{isomorphism classes of\ }
d\text {-dimensional modules with top\ } T\},$$  
sending $C$ to the class of $P/C$.  Clearly, the fibres of this map 
coincide with the orbits of the natural $\autlap$-action on
$\grasstd$.  

Recall moreover that by the radical layering of a module $M$ we mean the
sequence
$$\SS(M) = (M/JM, JM/J^2M, \dots, J^{L-1}M/J^L M, J^L M).$$
  Since we
identify semisimple modules with their isomorphism classes, we are, in
effect, dealing with a matrix of discrete invariants of $M$ keeping
count of the multiplicities of the simple modules in the various
semisimple layers $J^l M/J^{l+1}M$ of $M$.  Correspondingly, we consider
the following action-stable locally closed subvarieties of
$\grasstd$, which clearly cover the latter variety, namely  
$$\grassSS  = \{ C \in \grasstd \mid \SS(P/C) = \SS\}$$ 
for any
$d$-{\it dimensional semisimple sequence $\SS$ with top\/}
$T$; by this we mean any sequence  
$(\SS_0, \dots, \SS_L)$ of semisimple modules such that
$\SS_0 = T$ with each $\SS_l$ embedding into $J^lP/J^{l+1}P$, and
$\sum_{0 \le l \le L} \dim \SS_l = d$. 
\smallskip 

{\it Graded Grassmannians of modules with top $T$.\/} \ In keeping with
our goal of classifying graded modules, we restrict our focus to a
closed subvariety of $\grasstd$ paired with an action of the group of
homogeneous automorphisms of $P$ as follows:  Let
$M$ be a
$d$-dimensional graded module generated in degree zero with top
$T$, the latter being homogeneous of degree zero.  Clearly
$M$ has a graded projective cover $\pi: P \rightarrow M$, where $P =
\bigoplus_l P_l$ is equipped with the natural grading
and $\pi$ is homogeneous.  Thus we obtain a graded projective module
$P$ generated in degree zero such that, up to isomorphism,
$M$ equals
$P/C$, where $C$ is a homogeneous submodule of $P$.  If $M' = P/C'$
is another such module, then $M$ and
$M'$ are graded isomorphic if and only if there exists a homogeneous
automorphism $f: P
\rightarrow P$ with the property that $f(C) = C'$.  We
glean that classifying the $d$-dimensional graded modules with 
radical layering $\SS$ which are generated in degree zero (up to
graded-isomorphism) boils down to classifying the homogeneous points
$C \in \grassSS$.

To this end, we define 
$\ggrassSS$ to be the set of all those points $C \in \grassSS$ which
are homogeneous submodules of $JP$.  To see that
$\ggrassSS$ is projective, let $d_{li}$ be the multiplicity of
$S_i$ in $\SS_l$ and $P_{li}$ the $K$-subspace of $P_l$ generated
by all $p z_r$ with $1 \le r \le \frakt$ such that $p$ is a path (of length
$l$) ending in $e_i$.  Moreover, we denote the classical Grassmannian
of all $\bigl(\dim P_{li} - d_{li}\bigr)$-dimensional subspaces of
$P_{li}$ by
${\ggrassSS}_{li}$.  Then $e_i (P_l/C_l) \cong e_i \SS_l$ for $C$ as
above, whence
$$\ggrassSS =\{C\in \grassSS \mid C =
\bigoplus_{1 \le l \le L, 1 \le i \le n} C_{li} \text{\ with\ } C_{li}
\in {\ggrassSS}_{li}\}.$$
This subvariety of $\grassSS$ is isomorphic to a closed subvariety of
the product $\prod_{l,i} {\ggrassSS}_{li}$, and hence it is projective and
closed in $\grasstd$.  If, analogously, we define
$\ggrasstd$ to be the set of all points in $\grasstd$ which are
homogeneous submodules of
$JP$, then  
$$\ggrasstd  = \bigcup_{\SS} \ggrassSS$$
is a finite union of closed subvarieties of $\grasstd$, and 
consequently $\ggrasstd$ is projective as well.  

By $\gautlap$ we denote the subgroup of $\autlap$ consisting of
the homogeneous automorphisms of $P$.  Clearly, $\gautlap$ acts
morphically on the $\ggrassSS$ and $\ggrasstd$, and by the above
comments, the orbits are in one-to-one correspondence with the
isomorphism classes of graded $d$-dimensional modules (generated in
degree zero) having radical layering $\SS$ or top $T$, respectively.  
Note that
$\gautlap$ $\cong$
$\aut_\la (T)$ $\cong$ $\prod_{1 \le i \le n} \GL_{t_i}$.  Moreover,
we note that, given an isomorphism
$\overline{f}: P/C \rightarrow P/C'$ induced by
some $f \in \gautlap$, the ``distinguished" sequence $z_r + C$
of top elements of $P/C$ is mapped to a sequence of top elements of
$P/C'$, the
$f(z_r) + C'$ being again homogeneous
of degree zero.

Before we proceed, we give a pedantic reformulation of the
moduli problem tackled in Sections 3--5.   

\definition{Moduli problem for $d$-dimensional graded modules
generated in degree zero, with fixed top $T$}  For any point $C \in
\ggrasstd$, the module
$M = P/C$ has a natural grading, namely $M = \bigoplus_{0 \le l \le L}
M_l$, where
$M_l$ is the canonical image of $P_l$ modulo $C_l$; in particular, $M_0$
is the canonical image of $\bigoplus_{1 \le r \le \frakt} K z_r$, $M_1$
the image of $\bigoplus_{1 \le r \le \frakt} \sum_{\alpha \text{\ arrow}} K
\alpha z_r$, etc.  Clearly,
$M$ is generated in degree zero under this grading.  As we saw,
classifying the $d$-dimensional graded modules with top $T$
(resp\., with radical layering $\SS$) which are generated in degree zero
amounts to classifying the modules of the form $P/C$ with $C \in
\ggrasstd$ (resp\., with $C \in \ggrassSS$) up to graded-isomorphism.
  
Recall that a family of $d$-dimensional graded modules with
fixed top $T$ generated in degree zero is a family
$(\Delta, \delta)$ of modules each fibre of which has the listed
properties.  Two such families,
$(\Delta^1, \delta^1)$ and $(\Delta^2, \delta^2)$, parametrized by the
same variety
$X$, are equivalent if each $\Delta^1_x$ is graded-isomorphic to
$\Delta^2_x$.          

Our primary questions are:  When is there a
family that is universal for the families of equivalence classes of
graded $d$-dimensional modules generated in degree zero with top $T$; in
other words, when does our classification problem admit a fine moduli
space?  When does a coarse moduli space exist?  In case of existence,
what can be said about the geometry of such moduli spaces?
\enddefinition

\head 3. The distinguished affine cover of $\ggrasstd$ and a
pivotal family of graded modules
\endhead

The concepts introduced and discussed in this section will only be
relevant to proofs and examples; they will not appear in the statements
of the main results.

The affine cover of $\ggrasstd$ which we will present in the sequel
will provide the basis for the construction of a family of graded
modules which will be instrumental in solving the problems posed at
the end of Section 2.  This family will coincide with the universal
family in case a fine moduli space exists.   Moreover, being computable
from quiver and relations of
$\la$, this cover is the principal resource for the analysis of
examples.  We will briefly recall some concepts and facts provided in
\cite{\grassIII}, which will serve as our point of departure.  

\definition{Definitions 3.1}
  
{\rm (1)} An ({\it abstract\/}) {\it $d$-dimensional skeleton with top
$T$} is any sequence $\sigma = (\sigma^{(1)}, \dots, \sigma^{(\frakt)})$ of
sets
$\sigma^{(r)} \subseteq KQ e(r) \setminus I$ --  repetitions allowed  -- 
with $\sum_{1 \le r \le t} |\sigma^{(r)}| = d$, such that each $\sigma^{(r)}$ is
a nonempty set of paths starting in the vertex $e(r)$ which is closed
under right subpaths; the latter means that $pq \in \sigma^{(r)}$ implies
$q \in \sigma^{(r)}$.  (In particular, this condition guarantees that each
$\sigma^{(r)}$ contains the path $e(r)$ of length zero.) 

 To distinguish  between occurrences of the same path
$p$ in the intersection of two sets
$\sigma^{(r)}$ and $\sigma^{(s)}$ with $r \ne s$, but $e(r) = e(s)$, we tag
the elements of $\sigma^{(r)}$ with the superscript $r$ and write
$\sigma^{(r)}$ as $\{p^{(r)} \mid p \in \sigma^{(r)} \}$.  This
notational device makes it unambiguous to treat
$\sigma$ as a disjoint union $\bigsqcup_{1 \le r \le t} \sigma^{(r)}$.  

It is natural to view a skeleton $\sigma$ as a disjoint union of $\frakt$
tree graphs, each represented by some
$\sigma^{(r)}$, having ``root" $e(r)$ and edges labeled by the arrows
making up the paths in $\sigma^{(r)}$. 
 
{\rm (2)} Given a $d$-dimensional skeleton $\sigma$ with top $T$, a {\it
$\sigma$-critical pair\/} is a pair
$(\alpha, p^{(r)})$ with $1 \le r \le \frakt$, where
$\alpha$ is an arrow and $p^{(r)}$ a path in $\sigma^{(r)}$ such that
$\alpha p^{(r)}$ is a path of length at most $L$ in $KQ$ which does not
belong to $\sigma^{(r)}$.  Moreover, for any such $\sigma$-critical pair,
$\sigma(\alpha, p^{(r)})$ denotes the set of all paths $q \in \sigma$ (i.e., $q
= q^{(s)} \in \sigma^{(s)}$ for some $s$) with $\len(q) \ge \len(\alpha
p^{(r)})$ and  $\text{end}(q^{(s)}) = \text{end}(\alpha)$.   

{\rm (3)}  Let $M \cong P/C$ with  $C \in \grasstd$ and let $\sigma$ be a
$d$-dimensional skeleton with top $T$;  by $\sigma^{(r)}_l$ we denote the
set of paths of length $l$ in $\sigma^{(r)}$.  We say that $M$ {\it
has  skeleton
$\sigma$ relative to the presentation $\bigoplus_{1 \le r \le \frakt} \la z_r$
of $P$\/} in case the following holds: For each
$0 \le l \le L$, the union
$\bigcup_{1 \le r \le \frakt} \sigma^{(r)}_l (z_r + C)$ induces a basis for
$J^l M / J^{l+1}M$.  We call $\sigma$ {\it a skeleton of $M$\/} if
$\sigma$ is a  skeleton of
$M$ relative to some presentation of $P$.

{\rm (4)} Finally, we set 
$$\grassS = \{ C \in \grasstd \mid \sigma \text{\ is a skeleton of\ } M
\text{\ relative to\ }   P =  \bigoplus_{1 \le r \le \frakt} \la z_r \}$$
and
$$\ggrassS = \grassS \cap \ggrasstd.$$
\enddefinition

Whereas the concept of a skeleton of 
$P/C$ relative to the presentation $P = \bigoplus_{1 \le 1 \le r} \la
z_r$ is tied to the 
sequence $(z_r)_{r \le \frakt}$ of top elements, the set of {\it all\/}
skeleta of $M$ is an isomorphism
invariant of $M$.  In fact an abstract
$d$-dimensional skeleton $\sigma$ with top $T$ is a skeleton of $M$ if and
only if
\roster
\item"($\bullet$)" there exists a
sequence of top elements $m_1, \dots, m_\frakt$ of $M$ such that, for each
$l$,
$$\bigcup_{1 \le r \le \frakt} \sigma^{(r)}_l m_r \text{\ \ \ induces a basis
for\ \ \ } J^l M / J^{l+1}M.$$ 
\endroster 
In particular, $M$ is represented by {\it
some\/} point in $\grassS$ if and only if condition ($\bullet$)
holds.  (Yet, not all points in $\grass{\SS(M)}$ representing $M$ will
belong to $\grassS$ in general.) 

  Clearly, each
$d$-dimensional $\la$-module $M$ with top $T$ has at least one skeleton
relative to our fixed presentation $P = \bigoplus_{1 \le r \le \frakt} \la
z_r$.  In particular,
$\grasstd$ is covered by the locally closed subvarieties
$\grassS$.  More precisely, if $\SS$ is a
$d$-dimensional semisimple sequence, then $\grassSS$ is the union of
those charts $\grassS$ which have non-empty intersection with
$\grassSS$.  In fact,
$\grassSS$ is the union of those subvarieties $\grassS$ corresponding
to the skeleta $\sigma$ {\it compatible with\/} $\SS$ $=$ $(\SS_0, \dots,
\SS_L)$ in the following sense:  For each $l \le L$
and $i \le n$, the number of paths in $\bigsqcup_{r \le \frakt}
\sigma^{(r)}_l$ that end in the vertex $e_i$ coincides with the
multiplicity of the simple module $S_i$ in $\SS_l$.      

That $C$ be a point in $\grassS$ obviously entails the existence of
unique scalars
$c_{\alpha p^{(r)}, q^{(s)}}$ with the property that
$$\alpha p^{(r)} (z_r + C) =  \sum_{q^{(s)} \in \sigma(\alpha, p^{(r)})}
c_{\alpha p^{(r)}, q^{(s)}} q^{(s)} (z_s + C),$$ 
whenever $(\alpha,
p^{(r)})$ is a $\sigma$-critical pair.  Conversely, the isomorphism type
of $M = P/C$ is completely determined by the family of these scalars. 
Thus we obtain a bijection 
$$\psi: \grassS \rightarrow \bbA^N,\ \ \ \ C \mapsto c = \bigl(c_{\alpha
p^{(r)}, q^{(s)}}\bigr)_{(\alpha, p^{(r)})\ \sigma\text{-critical},\
q^{(s)} \in \sigma(\alpha, p^{(r)})}$$   
where $N$ is the {\it disjoint\/} union of the $\sigma(\alpha, p^{(r)})$.  
The latter sets of paths are not a priori disjoint; in the
index set $N$ we force disjointness by indexing the elements of
$\sigma(\alpha, p^{(r)})$ by the pertinent critical pair $(\alpha,
p^{(r)})$.  The following result from
\cite{\grassIII} summarizes the properties of the cover
$\bigl(\grassS \bigr)_\sigma$ which will be relevant here.   

\proclaim{Known facts 3.2 concerning the ungraded setting}  For
every
$d$-di\-men\-sion\-al skeleton
$\sigma$ with top $T$, the set $\grassS$ is a locally closed affine
subvariety of
$\grasstd$ which is isomorphic to a closed subvariety of $\bbA^N$ by
way of the map $\psi: \grassS \rightarrow \bbA^N$
described above.  Moreover, given any $d$-dimensional semisimple
sequence $\SS$ with top
$T$, the varieties $\grassS$, where $\sigma$ traces the skeleta
compatible with $\SS$, form an {\rm open} affine cover of $\grassSS$. 
{\rm{(}}On the other hand, the $\grassS$ fail to be open in $\grasstd$
in general.{\rm{)}}
\qed \endproclaim

Next we adjust the statements in 3.2 to the graded scenario.  The fact
that the affine sets $\ggrassS$ form an {\it open\/} cover of
$\ggrasstd$ in this case already anticipates the firmer grip we have on
the graded setting.

\proclaim{Facts 3.3 concerning the graded setting}  For
every $d$-dimensional skel\-eton
$\sigma$ with top $T$, the set $\ggrassS$ is an open affine
subvariety of $\ggrasstd$ which is isomorphic to a closed subvariety of
$\bbA^N$ by way of the restriction of $\psi$ to $\ggrassS$.  In
particular, given any
$d$-dimensional semisimple sequence $\SS$ with top
$T$, the varieties $\ggrassS$, where $\sigma$ traces the skeleta
compatible with $\SS$, form an open affine cover of $\ggrassSS$.
\endproclaim

\demo{Proof} First we observe that $\ggrassS$ coincides with the
intersection $\ggrasstd \cap \Schu(\sigma)$, where $\Schu(\sigma)$ is
the following open Schubert cell of the classical Grassmannian
$\Gr(d', JP)$ in which $\ggrasstd$ is located:  Namely, $\Schu(\sigma)$
coincides of those subspaces $C \subseteq JP$ for which $P = C \oplus
\bigoplus_{r \le \frakt, p^{(r)} \in \sigma^{(r)}} p^{(r)} z_r$; indeed
invoke gradedness to verify that $\sigma$ is a skeleton of $P/C$
precisely when the residue classes $p^{(r)} (z_r + C)$ with $r \le \frakt$
and $p^{(r)} \in \sigma^{(r)}$ form a basis for $P/C$. Thus the
$\ggrassS$ are open in $\grasstd$ and a fortiori in $\ggrasstd$.

The remainder of the claim will follow from the preceding facts if we
can show that the image of $\ggrassS$ under $\psi$ is a closed
subvariety of $\psi(\grassS)$ in $\bbA^N$.  But this is clear, since a
point 
$$c = \bigl(c_{\alpha
p^{(r)}, q^{(s)}}\bigr)_{(\alpha, p^{(r)})\ \sigma\text{-critical},\
q^{(s)} \in \sigma(\alpha, p^{(r)})}$$ 
in $\psi(\grassS)$ belongs to the image of $\ggrassS$ if and only if 
$c_{\alpha p^{(r)}, q^{(s)}} = 0$ whenever $q^{(s)}$ is strictly longer
than $\alpha p^{(r)}$. \qed \enddemo
     
One of the assets of this particular affine cover lies in its
computability.  Polynomials for the image of the map
$\psi: \grassS \rightarrow \bbA^N$ can easily be obtained from quiver and
relations of $\la$ (as was shown in \cite{\grassI}); a computer program
is being established by the authors of this paper.  We will see in
Section 5 that an even less labor-intensive algorithm yields polynomials
for
$\ggrassS$.  We will identify $C$ with $c$
whenever convenient.

The open affine cover
$\bigl( \ggrassS \bigr)_\sigma$ of
$\ggrasstd$ allows us to tie the graded modules $P/C$ with
$C \in \ggrasstd$ into a family, which will turn out to be universal (in
the sense of Section 2) in many cases.

\proclaim{Proposition 3.4}  There exists a family $(\Gamma, \gamma)$ of 
$d$-dimensional graded modules with top $T$, parametrized by
$\ggrasstd$, such that the fibre above any point $C \in \ggrasstd$,
endowed with the $\la$-module structure induced by $\gamma$, is $P/C$
with the grading inherited from $P$. \endproclaim 

\demo{Proof} For each $d$-dimensional skeleton $\sigma$ with top $T$, 
consider the following $d$-dimensional subspace 
$$V_\sigma = \bigoplus_{r \le t, p^{(r)} \in \sigma^{(r)}} K p^{(r)} z_r$$    
of the $K$-space $P$.  Note that the elements $z_r = e(r) z_r$
belong to each $V_\sigma$.  The trivial vector bundle
$\pi_\sigma: \ggrassS \times V_\sigma \rightarrow \ggrassS$ of rank
$d$ clearly becomes a family of graded modules when paired with the
algebra homomorphism $\gamma_\sigma: \la \rightarrow \End(\ggrassS
\times V_\sigma)$ which assigns to $\lambda \in \la$ the bundle endomorphism
that sends any $(C,x)$ to the pair $(C, y)$, where $y$ is the unique
element in $V_\sigma$ with $\lambda(x + C) = y + C$ in $P/C$. 

To glue these trivial bundles together along the intersections of
the affine patches $\ggrass{\sigma}$ of $\ggrasstd$, let $\rho$ and $\sigma$ be
two abstract
$d$-dimensional skeleta with top $T$, and let $g_{\rho, \sigma}:
\ggrass{\rho} \cap \ggrassS \rightarrow \GL_d$ be the map which sends
any point $C$ in the intersection $\ggrass{\rho} \cap \ggrassS$ to the
isomorphism $V_{\rho} \rightarrow V_\sigma$ that maps any $x \in V_{\rho}$
to the unique element $x' \in V_\sigma$ for which $x + C = x' + C$ in
$P/C$; note that $g_{\rho, \sigma}(C)$ is a well-defined $K$-space
isomorphism and $g_{\rho, \sigma}$ a morphism of varieties.   Moreover, our
construction entails that the
$g_{\rho, \sigma}(C)$ preserve the distinguished elements $z_r + C$ and
that the maps
$g_{\rho, \sigma}$, where $\rho$ and
$\sigma$ run through the $d$-dimensional skeleta with top $T$, satisfy the
coclycle condition, $g_{\sigma, \tau} (C) \circ g_{\rho, \sigma} (C) =
g_{\rho, \tau} (C)$, for all choices of $C \in \ggrass{\rho} \cap
\ggrassS \cap \ggrass{\tau}$.  Hence we obtain a vector bundle
$\Gamma$ over $\ggrassS$, which coincides with the above trivial bundles
on the affine patches of our cover.  The homomorphisms
$\gamma_\sigma$ are compatible with the gluing maps $g_{\rho, \sigma}$ and thus
yield an algebra homomorphism $\gamma: \la \rightarrow \End(\Gamma)$
which induces the $\la$-structure of $P/C$ on the fibre above any
point $C$ and, in particular, makes the $z_r(C)$ homogeneous of degree
zero.  \qed
\enddemo

We follow \cite{\New\rm, p. 37} in the concept of a {\it family
having the local universal property\/} relative to our moduli
problem.  The next observation serves merely as an
auxiliary to the proofs of our main results.    

\proclaim{Observation 3.5}  The family $(\Gamma, \gamma)$ of Proposition
{\rm 3.4} has the local universal property.  In other words, given an
arbitrary family
$(\Delta, \delta)$ of $d$-dimensional graded modules with top $T$ which
are generated in degree zero, any point
$x$ in the parametrizing variety $X$ of $\Delta$ has a neighborhood $U$
such that
$(\Delta|_U, \delta)$ is induced from $(\Gamma, \gamma)$ by way of some
{\rm(}not necessarily unique{\rm)} morphism $\tau: U \rightarrow
\ggrasstd$.  By the latter we mean that the family $(\Delta|_U, \delta)$
is equivalent to 
$\tau^*(\Gamma, \gamma)$. \endproclaim 

\demo{Proof}  Without
loss of generality $d \ge \dim T = \frakt$.  

Suppose $(\Delta, \delta)$ is 
any family of $d$-dimensional graded modules with top $T$ generated in
degree zero, $X$ being its parametrizing variety.  Let $x \in X$.  In
showing that, when restricted to a suitable neighborhood of $x$, this
family is induced by
$(\Gamma, \gamma)$, it is clearly harmless to
take $\Delta$ to be a trivial bundle $X \times K^d$.  We fix a sequence
$z_1, \dots, z_t$ of linearly independent elements of $K^d$ and assume
that, under the $\la$-structure induced by $\delta$ on any fibre
$\Delta_y$, these elements form a sequence $z_1(y), \dots,
z_\frakt(y)$ of top elements of $\Delta_y$; this assumption is permissible,
because
$(\Delta, \delta)$ is equivalent to a family $X \times K^d$ with this
property.   

Assume that the fibre of $\Delta$ above $x$ has skeleton
$\sigma$ relative to the sequence
$z_1(x), \dots, z_\frakt(x)$.  Then the set $Y$ of all
$y \in X$ such that $\Delta_y$ has skeleton $\sigma$ relative to
$z_1(y), \dots, z_\frakt(y)$ is an
open subset of $X$ containing $x$; indeed, the condition that the
elements $\delta(p^{(r)}) (z_r({y}))$ with $r \le \frakt$ and $p^{(r)}
\in \sigma^{(r)}$ be linearly independent is open.  We may thus
further assume that $X = Y$.  Thus, when verifying local
universality at $x$, we are essentially dealing with the trivial
subbundle $\ggrassS \times V_\sigma$ of
$(\Gamma, \gamma)$.  Identify the points
$C \in \ggrassS$ with the points $c \in \bbA^N$ described in the above
coordinatization of $\ggrassS$, and define a map
$\tau: X \rightarrow \ggrassS$ as follows:  For any $y \in X$, there
is a unique point $C = c(y) =  \bigl( c_{\alpha p^{(r)}, q^{(s)}} \bigr)
\in \ggrassS$ with the property
$$\alpha p^{(r)} z_r(y) =  
\sum_{q^{(s)} \in \sigma(\alpha, p^{(r)})} c_{\alpha p^{(r)}, q^{(s)}}
q^{(s)} z_s(y)$$ 
for all $\sigma$-critical pairs $(\alpha, p^{(r)})$.   
Setting $\tau(y) = c(y)$ yields a morphism $X \rightarrow
\ggrassS$, since, for each $r \le \frakt$, the map $\delta(-, z_r):  X
\rightarrow
\Delta = X \times K^d$, $y \mapsto \delta(y,z_r)$, is a morphism.  Thus
$\tau$ meets our requirement.  \qed
\enddemo

Finally, we record an immediate consequence of the preceding
observations which holds independent interest.

\proclaim{Corollary 3.6}  Each of the subvarieties $\ggrassSS$ of
$\ggrasstd$ is open and closed in $\ggrasstd$ and thus a union of
connected components. \endproclaim

\demo{Proof}  In Section 2, under the heading {\it Graded
Grassmannians\/}, we saw that, given any $d$-dimensional semisimple
sequence $\SS$ with top $T$, the set
$\grassSS$ is closed in $\ggrasstd$.  On the other hand, by Facts
3.3 above, $\grassSS$ is open, being a union of suitable charts
$\ggrassS$. \qed \enddemo

\head 4.  Sharpened version of Theorem A and a special case of
Theorem B \endhead

In light of Observation 3.5, a well-known criterion for the existence
of a coarse moduli space applies to the question of when the graded
$\la$-modules with fixed dimension and fixed top possess such a
(weakly) universal parametrizing space.  

  Recall that a
{\it categorical quotient\/} of $\ggrasstd$ by the action of $\gautlap$
is a variety
$\ggrasstd//\gautlap$ together with a morphism 
$$h: \ggrasstd
\rightarrow \ggrasstd//\gautlap$$ 
which is constant on the
$\gautlap$-orbits and has the property that every morphism defined on
$\ggrasstd$ and constant on these orbits factors uniquely through $h$.

\proclaim{Criterion 4.1} {\rm{(}}See \cite{\New\rm, Proposition
2.13}{\rm{)}}  The
$d$-dimensional graded modules with top $T$ {\rm{(}}resp., with radical
layering $\SS${\rm {)}} generated in degree zero possess a coarse moduli
space precisely when $\ggrasstd$ {\rm {(}}resp\.,
$\ggrassSS${\rm{)}} has a categorical quotient modulo
$\gautlap$ that separates $\gautlap$-orbits.   In case of existence,
this quotient {\it is\/} the coarse moduli space. 

In particular, closedness of the $\gautlap$-orbits in $\ggrasstd$ is a
necessary condition for the existence of a coarse moduli space.
\endproclaim

\demo{Proof}  For the final assertion, invest the fact that $\ggrasstd$
and all the $\ggrassSS$ are closed in $\grasstd$. \qed \enddemo
  
This citerion leads to the first of our main results, a
more detailed version of Theorem A.  As pointed out earlier,
each local graded module is graded-isomorphic to one generated in degree
zero. Therefore, given any simple $T$, the factor modules $P/C$ with $C
\in
\ggrasstd$, where $P$ is endowed with the natural grading, are
representative of all graded modules with top
$T$, up to graded-isomorphism.

\proclaim{Theorem 4.2}  Suppose that $\la = KQ/I$ is a finite
dimensional algebra which is graded by path lengths, and let $T$ be a
simple $\la$-module.  Then the projective variety $\ggrasstd$ is a fine
moduli space for the graded
$d$-dimensional $\la$-modules with top $T$, classifying them up to
graded-isomorphism.  Moreover, the family $(\Gamma, \gamma)$ constructed
in Proposition {\rm3.4} is universal.

A fortiori, given any $d$-dimensional semisimple sequence $\SS$ with
top $T$, the variety $\ggrassSS$ is a fine moduli space for
the graded $\la$-modules with radical layering $\SS$.

In both instances, the moduli space is a projective variety. 
\endproclaim

\demo{Proof}  Since the only
graded automorphisms of $P$ are multipications by nonzero scalars, 
$\gautlap \cong K^*$.  In particular, all homogeneous submodules of $P$
are invariant under $\gautlap$, which makes the $\gautlap$-orbits of
$\ggrasstd$ singletons.  Consequently, $\ggrasstd$ is its own
geometric (and a fortiori, categorical) quotient modulo $\gautlap$. 
Thus Criterion 4.1 shows $\ggrasstd$ to be a coarse moduli
space for the families of graded $d$-dimensional modules with
top $T$.

To see that $\ggrasstd$ is even a fine moduli space, we will check that
the family $(\Gamma, \gamma)$ of Proposition 3.4 is universal.  Let
$(\Delta,
\delta)$ be any family of graded modules of the form $P/C$ with $C \in
\ggrasstd$, parametrized by a variety
$X$ say;  again, the projective module $P$ is endowed with its natural
grading.  Define a map
$\tau: X
\rightarrow
\ggrasstd$ by sending any point $x \in X$ to the unique point $C \in
\ggrasstd$ having the property that the fibre of $\Delta$ above $x$ is
graded-isomorphic to $P/C$.  Due to the fact that
$\ggrasstd$ is already known to be a coarse moduli space for our
problem, the bijection $\alpha$ from the set of graded-isomorphism
classes of graded $d$-dimensional modules with top $T$ to the variety
$\ggrasstd$, defined by $[P/C] \mapsto C$, satisfies the
conditions of 1.6$'$ in \cite{\New\rm, p\. 24}.  The first of these
conditions guarantees that $\tau$ is a morphism of
varieties.  That $(\Delta,
\delta) \sim \tau^*(\Gamma, \gamma)$ is clear, as is uniqueness
of $\tau$ with this property.  This proves the first claim. 

Clearly, the restriction to the closed subvariety
$\ggrassSS$ of the above universal family parametrized by
$\ggrasstd$ is universal for the graded modules with radical
layering $\SS$.

That both $\ggrasstd$ and $\ggrassSS$ are projective was shown in
Section 2.   
\qed \enddemo

As pointed out before, given $\la$, and a semisimple sequence
$\SS$ with simple top $T$, the {\it full\/} class of modules with
radical layering $\SS$ need not have a moduli
space.  When it does, that moduli space need not be projective,  In
fact, any affine variety can be realized as a fine moduli space of the
form $\grassSS$; see \cite{\GeomI\rm, Theorem G} and \cite{\GeomIV\rm,
Corollary B}. 

While in the local case, that is, in the case of a simple top $T$,
classifying arbitrary graded modules is clearly equivalent to
classifying those that are generated in degree zero because every
local graded module is generated in a single degree, this is a
priori no longer true for graded modules with fixed, but unrestricted,
top.  As announced in the introduction, Section 6 will fill in the
gap.  There we will show that existence of a moduli space for
$d$-dimensional graded modules with top $T$ forces all of these modules
to be direct sums of local submodules.  This will reduce the general 
situation to the one addressed in the next theorem.

\proclaim{Theorem 4.3}  Let $\la = KQ/I$ be as in Theorem {\rm4.2}, and
$T$ any finite dimensional semisimple $\la$-module.  Moreover, let $\SS$
be any $d$-dimensional semisimple sequence with top $T$.  In the
following statements, the term ``graded module" will stand for ``graded
module generated in degree zero".
\smallskip

{\rm (1)}  If the $d$-dimensional graded $\la$-modules with radical
layering
$\SS$ have a coarse moduli space classifying them up to
graded-isomorphism, then all such modules are {\rm{(}}as graded objects{\rm{)}}
direct sums of local modules.  

A fortiori:  If the $d$-dimensional graded $\la$-modules with top
$T$ possess a coarse moduli space, they are all direct sums of graded local
components. 
\smallskip

{\rm (2)}  Conversely, suppose that all $d$-dimensional graded
$\la$-modules with radical layering $\SS$ {\rm{(}}resp., with top
$T${\rm{)}} are direct sums of graded local submodules.  Then there exists
a finite partition of the considered class of modules such that
the pertinent disjoint subclasses have fine moduli spaces providing
classification up to graded-isomorphism.
\endproclaim

\demo{Proof}  (1)  Suppose the graded $d$-dimensional modules with
radical layering
$\SS$ have a coarse moduli space.  Criterion 4.1 then forces all
$\gautlap$-orbits of the points $C \in \ggrassSS$ to be closed in
$\ggrasstd$ and hence in $\grasstd$ (see Section 2).

So we only need to show that, whenever $P/C$ with $C \in \ggrassSS$ is
not a direct sum of graded local modules, the orbit $\gautlap.C$ fails to be
closed in $\grasstd$.  Assume that $P/C = M \oplus N$, where
$M$, $N$ are graded and $M$ is indecomposable and nonlocal.  Without loss of
generality 
$C = U
\oplus V$ with $U \subseteq \bigoplus_{1 \le r \le u} \la z_r$ and $V
\subseteq \bigoplus_{u+ 1 \le r \le \frakt} \la z_r$ such that $M =
 (\bigoplus_{1 \le r \le u} \la z_r)/ U$ and $N = 
(\bigoplus_{u+1 \le r \le \frakt} \la z_r)/V$. In particular, $u > 1$. 
Given any $\tau
\in K^*$, we consider the automorphism $f_\tau \in \gautlap$ defined by
$f_\tau (z_1) =
\tau z_1$ and $f(z_r) = z_r$ for $r \ge 2$.  The curve $K^* \rightarrow
\gautlap.C$, given by $\tau \mapsto f_\tau (C)$, has a unique extension
$\PP^1 \rightarrow \overline{\gautlap.C}$ due to completeness of
$\overline{\gautlap.C}$.  We denote the value of this extension at
infinity by $C' =
\lim_{\tau \rightarrow \infty} f_\tau (C)$.  First we note that $C'$ is
again a $d'$-dimensional subspace of $JP$.  To see that $C'$ does not
belong to $\gautlap.C$, we let
$\pi_1: P \rightarrow \la z_1$ be the projection along
$\bigoplus_{r \ge 2} \la z_r$.  Setting $\mu = \dim \pi_1(C)$, we
pick elements
$b_1, \dots, b_\mu \in U$ such that $\pi_1(b_1), \dots, \pi_1(b_\mu)$
form a basis for $\pi_1(C) = \pi_1(U)$.  We supplement it with a basis
$b_{\mu +1},
\dots, b_{\nu}$ for $U \cap \Ker(\pi_1)$ to obtain a basis $b_1, \dots,
b_\nu$ for $U$.  Finally, we add on a basis $b_{\nu +1}, \dots, b_{d'}$
for $V$, which results in a basis 
$b_1, \dots, b_{d'}$ for $C$.  Clearly, $b_{\mu +1}, \dots, b_{d'}$ are
fixed by all
$f_\tau$, and hence belong to $C'$.  Moreover, the following
spaces are contained in $C'$ (cf\. \cite{\grassII\rm, Lemma 4.7}):  Namely,
the one-dimensional subspaces $\lim_{\tau \rightarrow \infty} f_\tau(K
b_r)$ for $r$ between $1$ and $\mu$.  If $b_r = \sum_{1 \le i \le \frakt}
\lambda_{ri} z_i$ with $\lambda_{ri} \in \la$, then the latter space
equals 
$$\lim_{\tau \rightarrow \infty} K \bigl(\lambda_{r1} z_1
+  \sum_{i \ge 2} (1/\tau) \lambda_{ri} z_i \bigr) = K \lambda_{r1}
z_1 =  K \pi_1(b_r).$$ 
Since the elements $\pi_1(b_1), \dots, \pi_1(b_\mu), b_{\mu +1}, \dots,
b_{d'}$ of $C'$ are linearly independent by construction, they
form a basis for $C'$.  Consequently, $C' = C_1 \oplus C_2
\oplus C_3$, where $C_1 = \pi_1(C) = \pi_1(U)$, $C_2 = U
\cap (\bigoplus_{2
\le r \le u} \la z_r) = \sum_{\mu + 1 \le r \le \nu} K b_r$, and
$C_3 = V$.  We thus obtain the following decomposition of the
$\la$-module
$P/C'$:
$$P/C' = \bigl(\la z_1/ C_1 \bigr) \oplus \biggl( \bigl(\bigoplus_{2 \le
r \le u} \la z_r \bigr) / C_2 \biggr) \oplus \biggl(\bigl(\bigoplus_{u+1
\le r\le d'} \la z_r \bigr)/ C_3 \biggr).$$
Clearly, all three summands are non-trivial, and the third equals
$N$.  This shows that the number of indecomposable summands of $P/C'$
exceeds that of $P/C$, whence $P/C
\not\cong P/C'$ as required, and the proof of (1) is complete.

For (2), we assume that all graded $d$-dimensional modules with top
$T$ are direct sums of locals, and consider the following equivalence
relation on the set of all partitions $(d_1,
\dots, d_\frakt)$ of $d$ with the property that $d_r
\ge 1$ for all $r$.  Namely, we call partitions $(d_1, \dots, d_\frakt)$
and $(d'_1, \dots, d'_\frakt)$ equivalent if, for every $a \in \NN$ and $i
\le n$, the number of $r \in \{1, \dots, \frakt\}$ with $e(r) = e_i$ and
$d_r = a$ equals the number of those $r$ for which $e(r) = i$ and $d'_r
= a$.  In the following, we will identify partitions of the described
type with their equivalence classes.  For each partition $(d_1, \dots,
d_\frakt)$, we thus obtain a class $\C(d_1, \dots, d_\frakt)$ of modules which
are direct sums $M_1 \oplus \cdots \oplus M_\frakt$, where $M_r$ is a
graded local module with top $S(r) = \la e(r)/J e(r)$.  In view of our
equivalence, the class of all graded
$d$-dimensional modules is the disjoint union of the classes $\C(d_1,
\dots, d_\frakt)$, where
$(d_1, \dots, d_\frakt)$ runs through the permissible partitions of
$d$.  By Theorem 4.2, the variety $\prod_{1 \le r \le \frakt}
{\operatorname{Gr-\frak{Grass}}}^{S(r)}_{d_r}$ is a fine moduli space
for the graded isomorphism classes of the objects in $\C(d_1, \dots,
d_\frakt)$.  The more restricted situation, where only the $d$-dimensional
graded modules with radical layering $\SS$ are assumed to be direct
sums of local submodules, is dealt with analogously.  This proves (2).
\qed
\enddemo     

\head 5.  Which varieties occur as fine moduli spaces of graded
modules?  Examples
\endhead

{\it Throughout this section, we let
$T$ be a simple $\la$-module: $T = S_0$ and $P = \la e_0$.  As before, we
assume $\la = KQ/I$ to be graded by path lengths.\/}

As we will see, the full spectrum of
possibilities already occurs in the situation of local
graded modules.  Namely, each irreducible projective variety
arises as a fine moduli space $\ggrassSS$, where $\SS$ is a
semisimple sequence with simple top $T$.  If the dimension of $\SS$ is
$d$, then $\ggrassSS$ is an irreducible component of $\ggrasstd$ in
this situation, showing that arbitrary irreducible projective varieties
can be realized as irreducible components of fine moduli spaces
$\ggrasstd$ (keep in mind that $\ggrassSS$ is always a union of
irreducible components of $\ggrasstd$).  We will use a family of
examples constructed by Hille in
\cite{\Hil} which, in turn, is based on a construction in
\cite{\GeomI\rm, proof of Theorem G}. Our arguments will illustrate the
general computational method sketched below. 

Our first set of examples is completely straightforward and does not
require any preparation. 

\proclaim{Examples 5.1}  Let $\SS$ be a semisimple
sequence with simple top $T$.  If $\SS_l = 0$ for $l \ge 2$,
then $\ggrassSS$ is either empty or a direct product of classical
Grassmannians $\Gr(u, K^v)$.
{\rm {(}}Note that in Loewy length two, $\ggrassSS = \grassSS$.{\rm {)}}

Conversely, every direct product of classical
Grassmannians occurs as a fine moduli space of local modules of
Loewy length two with fixed radical layering.
\endproclaim

\demo{Proof}  By our blanket hypothesis, $\SS_0 = T = S_1$.  For the
first statement, suppose
$\SS_1 =
\bigoplus_{1 \le i \le n} S_i^{m_i}$ and assume that $\ggrassSS \ne
\varnothing$.  Moreover, for $1\le i \le n$, let $\alpha_{i1}, \dots,
\alpha_{in_i}$ be the distinct arrows from the vertex $1$
to the vertex $i$ of $Q$.  Then $m_i \le n_i$, and
any $C$ in $\ggrassSS$ is a direct sum of subspaces $C_i$ of
$\bigoplus_{1 \le j
\le n_i} K \alpha_{ij} = e_iJe_1/e_iJ^2 e_1)$ of dimension $n_i - m_i$,
respectively.  Conversely, any such direct sum
of subspaces yields a point in $\ggrassSS$.  This
shows $\ggrassSS$ to be isomorphic to the direct
product $\prod_{1 \le i \le n} \Gr(n_i - m_i,
K^{n_i})$.

In light of the preceding paragraph, it is clear how to choose $Q$, so
as to realize any given product of Grassmannians $\Gr(d_i, K^{n_i})$.
\qed
\enddemo

Suppose $\SS$ is any $d$-dimensional semisimple sequence with top $T$
and $\sigma$ a skeleton compatible with $\SS$.  In \cite{\grassI}, we
described a method for determining the distinguished affine cover of
$\grassSS$  --  see Section 3  --  in the ungraded case.  The
algorithm provided there only requires a minor adjustment to yield the
corresponding cover for the graded version $\ggrassSS$; in fact,
computationally, the graded variant amounts to a significant reduction
of labor.   We will present the procedure for obtaining
polynomials defining
$\ggrassS$ in affine coordinates without proof, since the argument for
\cite{\grassI\rm, Theorem 3.14} is readily adapted to
the graded situation.   Returning to the notation of Section 3,
we first simplify the notation of a skeleton $\sigma$ with top $T =
S_1$, to reflect the fact that $\frakt = t_1 = 1$ in this section.  This
means $\sigma = \sigma^{(1)}$ and allows us to drop the superscript.  So, in
the present case, $\sigma$ is simply a set of paths in $KQ e_1$ of
cardinality $d$, which is closed under right subpaths.  As spelled out
in the proof of Facts 3.3, in dealing with graded
modules, we replace the set of paths
$\sigma(\alpha, p)$ for any $\sigma$-critical pair $(\alpha,p)$ by the following
subset $\operatorname{Gr-}\sigma(\alpha,p) \subseteq \sigma(\alpha,p)$. 
Namely, 
$$\operatorname{Gr-}\sigma(\alpha,p) = \{q \in \sigma(\alpha,p) \mid
\len (q) = \len (\alpha p)\}.$$ 
Then, clearly, any module $P/C$ with $C \in \ggrassS$ satisfies
$$\alpha p (e_1 + C) = \sum_{q \in
\operatorname{Gr-}\sigma(\alpha,p)} c_{\alpha p, q} q (e_1 + C)$$
for unique scalars $c_{\alpha p, q}$.
Modifying the notation of Section 3, we let $N$ be the disjoint union
of the sets $\operatorname{Gr-}\sigma(\alpha,p)$, where $(\alpha,p)$
traces the $\sigma$-critical pairs.  (Again, a priori, this union may fail
to be disjoint; we make it disjoint through suitable labeling.)  Then
the map
$$\psi: \ggrassS \rightarrow \bbA^N, \quad C \mapsto c =
\bigl(c_{\alpha p, q} \bigr)_{(\alpha,p)\ \sigma\text{-critical}, \ q \in
\operatorname{Gr-}\sigma(\alpha,p)}$$
is an isomorphism of varieties.  Our goal is to determine
polynomials whose zero locus in
$\bbA^N$ coincides with the image of $\psi$.  The polynomials we will
construct will be in variables $X_{\alpha p,q}$, where $(\alpha,p)$
traces the
$\sigma$-critical pairs and $q$ the corresponding sets
$\operatorname{Gr-}\sigma(\alpha,p)$.

\definition{5.2. The congruence relation induced by $\sigma$} Keeping $\sigma$
fixed, we consider the polynomial ring  
$$\A = \A(\sigma) := KQ[X_{\alpha p, q} \mid (\alpha, p)\
\sigma\text{-critical},\ q \in
\operatorname{Gr-}\sigma(\alpha,p)]$$ over the path algebra $KQ$.
On the ring $\A$, we consider congruence modulo the left ideal
$$\C = \C(\sigma) := \bigoplus_{2\le i\le n} \A e_i \ \ + \
\sum \Sb (\alpha,p)\ \sigma\text{-critical} \endSb \A \bigl( \alpha p \ -
\sum \Sb q\in \operatorname{Gr-}\sigma(\alpha,p) \endSb  X_{\alpha p, q} q
\bigr),$$ and denote this relation by $\seq$.
\enddefinition

In complete analogy with \cite{\grassI\rm, proof of Proposition 3.12}, one
shows that the quotient $\A/\C$ is a free left module
over  the commutative polynomial ring
$K[ X_{\alpha p, q}]$, the cosets $p + \C$ of the paths in $\sigma$ forming
a basis.  This means that, for any element $z \in \A$, there exist  
unique polynomials $\tau_q(X) = \tau^{z}_q(X)$ in 
$K[ X_{\alpha p, q}]$ such that 
$$z \seq \sum_{q \in \sigma} \tau_q(X)\ q.$$

We explain how to obtain the $\tau_q(X)$ in case $z = p$ is a path,
this being sufficient for dealing with arbitrary elements $z \in \A$. 
If $p$ starts in a vertex different from $e_1$, set
$\tau_q(X) = 0$ for all $q$.  Now suppose that $p = pe_1$, and let
$p_1$ be the longest right subpath of $p$ that belongs to $\sigma$; this
path may have length zero, that is, coincide with $e_1$.  If $p_1 = p$,
set $\tau_p(X) = 1$ and $\tau_q(X) = 0$ for $q \ne p$; otherwise, write
$p = p'
\alpha p_1$ for some $\sigma$-critical pair $(\alpha, p_1)$ and some left
subpath $p'$ of $p$, potentially of length zero.  Then $z \seq 
\sum_{q \in \operatorname{Gr-}\sigma(\alpha,p)} X_{\alpha p,q} p' q$. 
Iterate this step for each of the paths $p'q$ appearing in the
latter sum, noting that they all have strictly longer right subpaths
in $\sigma$ than does $p$, while having the same length as $p$.  Thus an
inductive procedure will take us to the desired normal form of $z$ under
$\seq$.  

\definition{5.3 Polynomials for $\ggrassS$} Let
$\R$ be any finite generating set for the left ideal $Ie_1 + \cdots
+Ie_t$ of
$KQ$  --  note that such a generating set always exists since all paths
of lengths $\ge L+1$ belong to $I$.  For each
$\rho\in \R$, let
$\tau^{\rho}_q(X)$, $q \in \sigma$, be the unique polynomials in $K[
X_{\alpha p, q}]$ with
$$\rho \  \seq \ \sum_{q\in \sigma} \tau^{\rho}_q(X)\ q,$$ 
as guaranteed by 5.2.  Then the zero locus
$V\bigl(\tau^{\rho}_q(X) \mid \rho \in \R,\, q \in \sigma \bigr)$ in
$\bbA^N$ of these polynomials is the image of $\ggrassS$
under the isomorphism  $\psi: \ggrassS \rightarrow \Img(\psi) \subseteq
\bbA^N$.  
\enddefinition

From Examples 5.1 we already know that, for $n \ge 0$, projective
$n$-space can be realized as a fine moduli space $\ggrassSS$ for a
suitable semisimple sequence $\SS$ with simple top.  In
preparation for the examples announced at the beginning of this
section, we again realize $\PP^n$ in the form $\ggrassSS$, but this
time in a more ``ample" setting that will allow us to modify
$\ggrassSS$ by means of additional relations factored out of the
pertinent path algebra $KQ$.  
 Since the above method allows for verification of our
claims in very elementary terms, we will include brief arguments for the
convenience of the reader.

\definition{Examples 5.4} \cite{\Hil}  Let $Q$ be the quiver 

$$\xymatrixrowsep{2.0pc}\xymatrixcolsep{6pc}
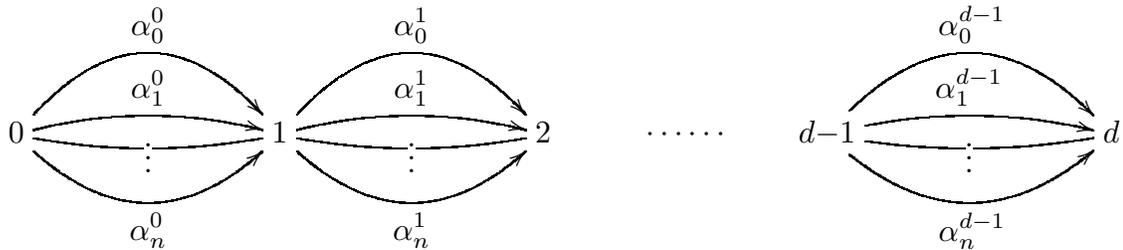
\xymatrix{
0 \ar@/^6ex/[r]^{\alpha^0_0} \ar@/^/[r]^{\alpha^0_1} \ar@{}@/_1ex/[r]|{\vdots}
\ar@/_5ex/[r]_{\alpha^0_n}  &1 \ar@/^6ex/[r]^{\alpha^1_0} \ar@/^/[r]^{\alpha^1_1}
\ar@{}@/_1ex/[r]|{\vdots}
\ar@/_5ex/[r]_{\alpha^1_n}  &2 \ar@{}[r]|{\displaystyle\cdots\cdots} &d{-}1
\ar@/^6ex/[r]^-{\alpha^{d-1}_0}
\ar@/^/[r]^{\alpha^{d-1}_1} \ar@{}@/_1ex/[r]|{\vdots}
\ar@/_5ex/[r]_-{\alpha^{d-1}_n}  &d
}$$

\noindent and $\la = KQ/I$, where $I$ is generated by all differences
$\alpha_i^{r} \alpha_j^{r-1} - \alpha_j^{r} \alpha_i^{r-1}$ for $0
\le i,j \le n$ and $1 \le r \le d-1$.  Consider the $(d+1)$-dimensional
semisimple sequence 
$\SS = (S_0, \dots, S_d)$.  Again, $P = \la e_0$. 
We will verify that $\ggrassSS \cong \PP^n$ for any choice of $d$.

If $d = 1$, the modules with radical layering $\SS$ are uniserial of
length $2$, and by 5.1, we obtain $\ggrassSS \cong \Gr(n, K^{n+1})
\cong \PP^n$ as desired.  The isomorphism $\Gr(n, K^{n+1}) \cong
\PP^n$ being non-canonical, it will be helpful to
specify a concrete incarnation, say $F: \PP^n
\rightarrow \ggrassSS$,  
$$(k_0: k_1: \cdots : k_n) \mapsto \sum_{0 \le i,j \le n} \la \bigl(k_i
\alpha_j^0 - k_j \alpha_i^0 \bigr).$$
Note that, for $k_0 = 1$, this latter submodule of $P$ equals $\sum_{1
\le i
\le n} \la \bigl(\alpha_i^0 - k_i \alpha_0^0 \bigr)$.   

Now let $d \ge 1$ be arbitrary.  We ascertain that every module
$M$ with radical layering $\SS$ is completely determined by
the factor module $M/J^2 M$ of length $2$ with radical layers
$(S_0, S_1)$.  Let $x$ be a top element of $M$.  Due to the symmetry of
our example, we may assume without loss of generality that $\alpha_0^0
x \ne 0$, meaning that
$M/J^2M  \cong  P/ F(k_0: k_1: 
\dots: k_n)$, with $k_0 = 1$ and suitable $k_i \in K$ for $i \ge 1$. 
Then it is readily checked
 --  for details consult the following paragraph  --  that
$M$ has skeleton
$\sigma = \{e_0, \alpha_0^0, \alpha_0^1 \alpha_0^0, \dots, p\}$, where
$p = \alpha_0^{d-1} \cdots \alpha_0^0$, and $M$ is completely determined
by the scalars $k_i$.  Indeed, $\alpha_{i_r}^{r} \cdots
\alpha_{i_0}^0 x = k_{i_r}
\cdots k_{i_0} \alpha_0^r  \cdots \alpha_0^0 x$ for any $r \le
d-1$ and any choice of $i_0, \dots, i_r \in \{0, \dots, n\}$. 
Conversely, it is clear that every module of length
$2$ with radical layers $(S_0, S_1)$ occurs as a quotient of a module
$M$ with radical layering $\SS$.  Thus, again, $\ggrassSS \cong
\ggrass{(S_0, S_1)} \cong \PP^n$, an isomorphism being available as in
the case $d = 1$. 

For justification, we display the congruences
determining $\grassS$, where $\sigma$ is the skeleton consisting
of all right subpaths of $p = \alpha_0^{d-1} \cdots \alpha_0^0$ as
above.  As a left ideal, $I$ is generated by the differences of the form
$\alpha_{i'_s}^{s} \alpha_{i'_{s-1}}^{s-1} \cdots \alpha_{i'_r}^{r} - 
\alpha_{i_s}^{s} \alpha_{i_{s-1}}^{s-1} \cdots \alpha_{i_r}^{r}$ for $0
\le r < s \le d-1$, where $(i'_s, \dots, i'_r)$ is any permutation of
$(i_s, \dots, i_r)$.  Moreover, we have (dn) $\sigma$-critical pairs, 
$\bigl(\alpha_i^{r}, \alpha_0^{r-1} \cdots \alpha_0^0
\bigr)$ for $0 \le r \le d-1$ and $0 \le i \le n$, giving rise to the 
basic congruences $\alpha_i^{r} \alpha_0^{r-1}
\cdots \alpha_0^{0} \seq X_i^{r} \alpha_0^{r} \alpha_0^{r-1} \cdots
\alpha_0^{0}$.  As a result of substituting them into any path of the
form $q =
\alpha_{i_r}^{r} \cdots \alpha_{i_0}^{0}$ of $\la$, we obtain the
following list of congruences:
$\alpha_{i_r}^{r} \cdots \alpha_{i_0}^{0}  \seq  X_{i_r}^r
\cdots X_{i_0}^0 \alpha_0^{r}  \cdots \alpha_0^{0}$.  The variables
$X_0^0, \dots, X_n^0$ can be chosen freely.  The relations in $I$ thus
yield $X_i^r = X_i^0$ for $0 \le i \le n$ and 
all $r$, in accordance with the previous paragraph.
\enddefinition

The next examples finally show that, indeed, every irreducible
projective variety takes on the form $\ggrassSS$ for some graded algebra
$\la$ and a suitable semisimple sequence $\SS$ with simple top.   They
will thus confirm the assertion we made at the beginning of the
section. 

\definition{Examples 5.5} \cite{\Hil}  Let $V \subseteq \PP^n$ be an
irreducible projective variety, determined by homogeneous polynomials
$f_1, \dots, f_s \in K[X_0, \dots, X_n]$ say.  Suppose moreover that
the degrees $d(1), \dots, d(s)$ of these polynomials are bounded from
above by
$d$.  To realize $V$ in the form $\ggrassSS$, let $\la = KQ/I$, where
$Q$ is the quiver of Examples 5.4.  The ideal $I$ of $KQ$ is generated
by the relations listed in 5.4, next to the
following additional ones, labeled $g_1, \dots, g_s$, one for every
polynomial $f_r$:  Write each monomial occurring in
$f_r$ in the form $X_{i_{d(r)}} \dots X_{i_1}$ with $i_j \in \{0,
\dots, n\}$, where the order of the factors is
irrelevant, and replace each variable $X_{i_j}$ by the
arrow $\alpha_{i_j}^j$. This process results in a homogeneous linear
combination $g_s$ of paths in $Q$. 

As before, $\SS = (S_0, \dots, S_d)$ and $P = \la e_0$.

In light of Examples 5.4, $\grassSS$ is a subvariety of $\PP^n$;
indeed, since our present ideal $I$ contains that of 5.4, each module
$M$ with radical layering
$\SS$ is uniquely determined by $M/J^2M$ with radical layering $(S_0,
S_1)$.  Moreover, in 5.4, we explicitly provide an assignment sending
$\overline{k} \in \PP^n$ to a module with radical layering $(S_0,
S_1)$.  

So we only need to show that the given variety $V$ coincides with the
set of those points
$\overline{k} = (k_0:
\dots : k_n) \in \PP^n$ for which there exists $C \in
\grassSS$  --  over the present incarnation of $\la$  --  with the
property that the factor $M/J^2M$ of  $M = P/C$ is determined by
$\overline{k}$.  Without loss of generality, we may assume that
$V$ is not contained in the hyperplane $k_0 = 0$, whence the
affine variety
$V_*$ obtained from $V$ by dehomogenizing at the variable $X_0$
completely determines $V$.  Hence our task is reduced to showing
that the set of those points $C \in \ggrassSS$ for which $\alpha_0^0
(e_0 + C) \ne 0$ in $P/C$ coincides with $V_*$.  The latter
is exactly the affine patch $\ggrassS$ where $\sigma$ is the skeleton
introduced in 5.4.  

The basic congruences leading to polynomials for $\ggrassS$ are as
listed in the last paragraph of 5.4, and since the relations in 5.4 are
among those we factored out of
$KQ$ in the present example, we again obtain 
$X_{i}^r = X_{i}^0$ for $0 \le i \le n$ and all 
$r \in \{0, \dots, d-1\}$.  Moreover, as before, any path $q =
\alpha_{i_r}^{r} \cdots \alpha_{i_0}^{0}$ in $\la$ is congruent to
$X_{i_r}^r \cdots X_{i_0}^0 \alpha_0^{r}  \cdots \alpha_0^{0}$.  Under
the legitimized identification of $X_{i}^r$ with $X_i$ for any $r$, 
substitution of the basic congruences into the relations $g_r$ thus
leads to the dehomogenizations of the polynomials $f_r$ relative to the
variable $X_0$. \qed
\enddefinition

\head 6.  Appendix:  Graded modules generated in mixed degrees \endhead

The purpose of this appendix is to show that, in Theorem 4.3, our
restriction to graded modules generated in degree zero is superfluous.

Again, we let $\la = KQ/I$ be graded by path lengths, and 
$$T =
\bigoplus_{1 \le i \le n} S_i^{t_i} = \bigoplus_{1
\le r \le \frakt} \la e(r)/ J e(r),$$ 
where $\frakt = \sum_i t_i$.  As before,
$P = \bigoplus_{1 \le r \le \frakt} \la z_r \rightarrow T$ is a projective
cover of $T$ sending the top elements $z_r = e(r) z_r$
to the residue classes $e(r) + J e(r)$, but now
we assume the simple modules $\la e(r) / J e(r)$ to be homogeneous
of degree $h(r)$, respectively, and call $h = \bigl(h(1), \dots, h(\frakt)
\bigr)$ the {\it degree vector\/} of $T$.  Correspondingly, we choose the
grading of
$P$ so as to make the above projective cover homogeneous; in other
words, we assume
$z_r$ to be homogeneous of degree
$h(r)$ for $r \le \frakt$.  Any factor module $M$ of $P$ by a homogeneous
submodule contained in $JP$ is then said to have {\it top degree vector\/} $h$. 
A full sequence of top elements of $M$ consists of a generating set
$m_1, \dots, m_{\frakt}$ such that $m_r = e(r) m_r$ is homogeneous of degree
$h(r)$; this provides the framework for carrying over the concept of a
skeleton of $M$.  Our goal is to verify that, in case the graded
$d$-dimensional
$\la$-modules with top
$T$ and top degree vector $h$ have a coarse moduli space, all of the
considered modules are direct sums of graded local submodules generated in the
degrees $h(1), \dots, h(\frakt)$.  We have already dealt with the special case $h
= (0, \dots, 0)$ in Theorem 4.3.  In the following, we will
outline the considerations required to adjust the argument to an
arbitrary top degree vector. 

To that end, we modify the definitions of Sections 2 and
3 as follows:   Let $\SS$ be a $d$-dimensional semisimple sequence
with top $T$ and, as in Section 2, denote by $d_{li}$ the
multiplicity of the simple module $S_i$ in $\SS_l$.  Moreover, let
$P_{li}$ be the $K$-subspace of $P$ generated by all elements of
the form $p z_r$, where $p$ is a path of length $l - h(r)$ ending
in the vertex $e_i$.  Again $\ggrassSS_{li}$ denotes the classical
Grassmannian of all $\bigl(\dim P_{li} - d_{li} \bigr)$-dimensional
subspaces of $P_{li}$.  But now, we define
$$\ggrassSS = \{C \in \grassSS \mid C = \bigoplus_{1 \le l \le L,\;
1 \le i \le n} C_{li} \text{\ with\ } C_{li} \in \ggrassSS_{li}
\},$$
and let $\ggrasstd$ be the union of the $\ggrassSS$, where $\SS$
runs through the $d$-dimensional semisimple sequences with top $T$.  
The same arguments as used in Section 2 guarantee projectivity of these
subvarieties of $\grasstd$.  The acting group $\gautlap$ is once more
the group of homogeneous $\la$-automorphisms of $P$; note, however, that
in mixed degrees, this group may have nontrivial unipotent radical.

Defining $\ggrassS$ as the intersection $\grassS \cap \ggrasstd$  -- 
this conforms with Section 3  --  one establishes analogues of Facts
3.3; the only adjustment required in the proof concerns the last
sentence:  Namely, in the present situation, the image of $\ggrassS$
under
$\psi$ consists of those points 
$$c = \bigl(c_{\alpha
p^{(r)}, q^{(s)}}\bigr)_{(\alpha, p^{(r)})\ \sigma\text{-critical},\
q^{(s)} \in \sigma(\alpha, p^{(r)})}$$
in $\psi(\grassS)$ for which $c_{\alpha p^{(r)}, q^{(s)}} = 0$
whenever 
$$\len(q^{(s)}) + h(s) \ne \len(\alpha p^{(r)}) + h(r).$$ 
Observation 3.4 remains unchanged; keep in mind that, for $C \in
\ggrasstd$, the grading of the factor module $P/C$ inherited
from $P$ now has top degree vector $h = \bigl(h(1), \dots, h(\frakt)
\bigr)$.  Observation 3.5 should be replaced by the remark that
the family of Proposition 3.4 has the local universal property
relative to families of graded $d$-dimensional modules with top
$T$ and the specified top degree vector $h$.  As in Criterion 4.1,
this setup allows us to apply \cite{\New\rm, Proposition 2.13} to
conclude that existence of a coarse moduli space for the
considered graded modules with top $T$ (resp., with radical layering
$\SS$) implies closedness of the $\gautlap$-orbits of
$\ggrasstd$ (resp., of $\ggrassSS$).   A replica of the argument
backing part (1) of Theorem 4.3 finally shows that the latter
closedness conditions force all $d$-dimensional graded modules
with top degree vector $h$ and top $T$ (resp., all $d$-dimensional
graded modules with top degree vector $h$ and radical layering
$\SS$) to be direct sums of graded local summands generated in the
degrees $h(1), \dots, h(\frakt)$.  But, as we already saw, this legitimizes
waiving of the hypothesis that
$h(r) = 0$ for all
$r$ in Theorem 4.3.  

{\it Conclusion\/}:  For any choice of $h$, the $d$-dimensional
modules  with top
$T$ and top degree vector $h$ have a coarse moduli space if and only if
the $d$-dimensional modules with top $T$ generated in degree
$0$ do, and in case of existence, the two moduli spaces
coincide.

\enddocument